\documentclass{amsproc}

\usepackage{a4wide}
\usepackage{amsmath}
\usepackage{enumerate}

\usepackage{graphicx, floatflt}
\usepackage[english]{babel}

\usepackage{amsmath,amsthm}
\usepackage{amsfonts}
\usepackage{amssymb}
\usepackage{longtable}
\usepackage[matrix,arrow,curve]{xy}

\newtheorem{theorem}{Theorem}

\newtheorem{prop}{Proposition}
\newtheorem{defin}{Definition}
\newtheorem{lem}{Lemma}

\newtheorem{ex}{Example}{\bf}{\rm}

{\bf}{\rm}{\rm}

\newtheorem*{theorem*}{Theorem}
\newtheorem*{cor*}{Corollary}

\def\Real{\mathbb{R}}

\def\V{\mathcal{V}}

\def\D{\mathcal{D}}

\def\F{\mathcal{F}}

\def\GLn{\text{\rm GL}(n,\Real)}
\def\On{\text{\rm O}(n)}
\def\SLn{\text{\rm SL}(n,\Real)}
\def\SLnn{\text{\rm SL}(n,\Real)\times\mathbb{Z}_2}

\def\GL1{\text{\rm GL}(1,\Real)}
\def\O1{\text{\rm O}(1)}

\def\Diff{\text{\rm Diff}}

\def\pt{\mathop\text{\rm pt}\nolimits}

\usepackage{stackrel}

\newcommand{\be}{\begin{equation}}
\newcommand{\ee}{\end{equation}}

\let\leq=\leqslant
\let\geq=\geqslant

\begin{document}

\title{On Losik classes of diffeomorphism pseudogroups}

\author{Yaroslav V. Bazaikin}\thanks{$^1$ Department of Mathematics, Faculty of Science, Jan Evangelista Purkyňe University in
	Ústí nad Labem, Pasteurova 3632/15, 400 96 Ústí na Labem,   Czech Republic}

\author{Yury D.  Efremenko}\thanks{$^2$Novosibirsk State University, Pirogova 1, 630090 Novosibirsk, Russia}

\author{Anton S. Galaev}\thanks{$^3$Department of Mathematics, Faculty of Science, University of Hradec Kr\'alov\'e,  Rokitansk\'eho 62, 500~03 Hradec Kr\'alov\'e,  Czech Republic\\
	Corresponding author e-mail: anton.galaev(at)uhk.cz}

\maketitle

\begin{abstract}

Let $P$ be a pseudogroup of local diffeomorphisms of an $n$-dimensional smooth manifold $M$. Following Losik we consider characteristic classes of the quotient $M/P$ as elements of the de~Rham cohomology of the second order frame bundles over $M/P$ coming from the generators of the Gelfand-Fuchs cohomology. We provide explicit expressions for the classes that we call Godbillon-Vey-Losik class and the first Chern-Losik class. Reducing the frame bundles we construct bundles over $M/P$ such that the Godbillon-Vey-Losik class is represented by a volume form on a space of dimension $2n+1$, and the first Chern-Losik class is represented by a symplectic form on a space of dimension $2n$. Examples in dimension 2 are considered.      

\vskip0.5cm


\vskip0.5cm




\end{abstract}

\section*{Introduction}

The present paper is a continuation of the works \cite{BGG,BG} developing Losik's approach to characteristic classes of foliations and diffeomorphism groups \cite{L90,L94,L10,L15}. For information on foliations,  diffeomorphism groups, and their characteristic classes see, e.g., \cite{Cr03,psgroups,Hurder00,Nar,NarYaz}  and the references therein. Let us mention another recent approach 
\cite{SZ1,SZ2}.   

A $D_n$-space $X$ introduced by Losik \cite{L94} is a set with a smooth structure generalizing the smooth structure of an $n$-dimensional manifold. Typically $X$ is the quotient $M/P$ of an $n$-dimensional manifold $M$ by a diffeomorphism pseudogroup $P$. An important example is the leaf space $N/\F$ of a regular foliation $\F$ on a manifold $N$, where $N/\F$ may be considered as the quotient of a complete transversal $M$ by the holonomy pseudogroup $P$. The characteristic classes 
of a quotient $X$ are  elements of the de~Rham cohomology of the second order frame bundles over $X$ coming from the generators of the Gelfand-Fuchs cohomology. These classes may be mapped to the usual characteristic classes of diffeomorphism groups or foliations, but then one looses information. For example, the Godbillon-Vey class of any Reeb  foliation $\F$ on the 3-dimensional sphere $S^3$ is trivial, while the corresponding characteristic class of the $D_1$-space $S^3/\F$ may be non-trivial~\cite{BGG}. 

For a quotient $M/P$ we consider first of all the Godbillon-Vey-Losik (GVL) class and the first Chern-Losik class (CL class). The GVL and CL classes are the cohomology classes given by a $(2n+1)$-form and a $2$-form on the spaces $S_2(M/P)/\On$ and $S_2(M/P)/\GLn$, respectively. Here  $S_2(M/P)$ is the bundle of frames of second order over $M/P$. The dimension of the spaces $S_2(M/P)/\On$ and $S_2(M/P)/\GLn$ are respectively $\frac{1}{2}n(n^2+2n+3)$ and $\frac{1}{2}(n^2+n+2)$. 
We find simplified versions of these classes and  consider the GVL class as the cohomology class of a volume form on a $D_{2n+1}$-space $A(M/P)$, also we consider the CL class as the cohomology class of a symplectic form on a $D_{2n}$-space $B(M/P)$. 

Let us give the explicit construction.
Let $U\subset M$ be a coordinate neighborhood with the coordinates $y^1,\dots, y^n$.   Let $A(U)=U\times \Real\times \Real^n$ be the bundle over $U$ with the coordinates  $y^1,\dots, y^n,x,y_1,\dots y_n$. 
Suppose that $V$ is another  coordinate neighborhood with the coordinates $\alpha^1,\dots, \alpha^n$.
We extend a  local diffeomorphism from $U$ to $V$ given by $\alpha^i=\alpha^i(y^1,\dots,y^n)$  to a local diffeomorphism from $A(U)$ to $A(V)$ by setting
$$\alpha^i=\alpha^i(y^1,\dots,y^n),\quad \beta=x+\ln\left|\det\left(\frac{\partial \alpha^i}{\partial y^j}\right)\right|,\quad \alpha_k=\frac{\partial y^p}{\partial \alpha^j }
\frac{\partial^2 \alpha^j }{\partial y^p \partial y^q}\frac{\partial y^q}{\partial \alpha^k}+
\frac{\partial y^j}{\partial \alpha^k }
y_{j}.$$ This defines a bundle $A(M)$ over $M$.   
For each $A(U)$, consider the differential form
$$-dx\wedge (dy_i\wedge dy^i)^n\in\Omega^{2n+1}(A(U)).$$ 
The above coordinate transformations transform these forms to each other. 
This gives a $(2n+1)$-form $gvl$ on $A(M)$. Each diffeomorphism of $M$ is extended to a diffeomorphism of $A(M)$ in the obvious way.
The form $gvl$ is invariant under all such diffeomorphisms, hence it defines a $(2n+1)$-form on $A(M/P)$, which we denote by the same symbol. We consider the cohomology class of $gvl$ in $H^{2n+1}(A(M/P))$ as a simplified version of the GVL class. This GVL class is trivial if and only if there exists a $P$-invariant $2n$-form $\omega$ on $A(M)$ such that $d\omega=gvl$.

As above, for each $U$ consider the bundle $B(U)=U\times \Real^n$ with the coordinates $y^1,\dots, y^n,y_1,\dots y_n$. We obtain the bundles $B(M)$ and $B(M/P)$.     The 2-form $$dy_i\wedge dy^i$$
on $B(U)$ defines a 2-form on $B(M/P)$. We consider the cohomology class of this form as a simplified version of the first CL class. 




We define characteristic classes of a vector field $V$ on $M$ as the characteristic classes of the corresponding pseudogroup $P$. In Section \ref{secex}, we consider examples of vector fields on the open disc $D^2\subset\Real^2$. Let $f$ be a smooth even function on $(-1,1)$, non-vanishing on $(0,1)$. Then  the function $f(r)$, where $r=\sqrt{(y^1)^2+(y^2)^2}$, is smooth on $D^2$. The first CL class of vector field $V_1$ with the components $(f(r)y^1,f(r)y^2)$ is non-trivial if and only if $f(0)=0$. This means that the triviality of the first CL class of $V_1$ is equivalent to the non-degeneracy of the singular point $(0,0)$. The first 
CL class of the vector field $V_2$ with the components $(f(r)y^1,-f(r)y^2)$ is trivial for each $f$; this is not a surprise, since the quotient $D^2/V_2$ is just the segment~$[0,1)$. Thus the characteristic classes under consideration provide information about singularities.

\section{Definition of characteristic classes}

Here we give a very short introduction to Losik's ideas. The details may be found in \cite{L90,L94,L10,L15}. 
Let $M$ be an $n$-dimensional smooth manifold. Denote by $S(M)$ the space of frames of infinite order over $M$, see \cite{BR} for details.
Consider also the reduced bundles $S(M)/\GLn$ and $S(M)/\On$. Let $P$ be a pseudogroup of local diffeomorphisms of $M$. Losik \cite{L94} introduced on $M/P$ a $D_n$-structure generalizing the notion  of the smooth structure of an $n$-dimensional smooth manifold. For simplicity we call $M/P$ just a quotient and denote it by $X$. A morphism of two quotients is a morphism of the corresponding pseudogroups defined by Haefliger \cite{psgroups}.    
The action of $P$ on $M$ is extended to an action of $P$ on the bundles  $S(M)$, $S(M)/\GLn$ and $S(M)/\On$. Let $S(M/P)=S(M)/P$ and define the spaces 
$S(M/P)/\GLn$ and $S(M/P)/\On$ in a similar way. A differential form on $S(M/P)$ is a $P$-invariant differential form on $S(M)$.
Let $\Diff(\Real^n)$ be the pseudogroup of all local diffeomorphisms of $\Real^n$. Let $\pt_n=\Real^n/\Diff(\Real^n)$.
There is a unique morphism $$X=M/P\to\pt_n$$ inducing the  morphisms
$$S(X)\to S(\pt_n),\quad S(X)/\GLn \to S(\pt_n)/\GLn,\quad S(X)/\On\to S(\pt_n)/\On.$$
There are the following isomorphism of the de~Rham cohomology of the bundles over $\pt_n$ and  Gelfand-Fuchs cohomology of the Lie algebra $W_n$ of formal vector fields on $\Real^n$:
\begin{equation}\label{isomH}
H^*(S(\pt_n))/G\cong H^*(W_n,G),
\end{equation} 
where $G$ is a trivial group or one of $\GLn$, $\On$.
This gives the characteristic morphisms
\begin{align}\label{chi}
H^*(W_n)&\to H^*(S(X)),\\ \label{chi'} H^*(W_n,\GLn)&\to
H^*(S(X)/\GLn),\\ \label{chi''} H^*(W_n,\On)&\to
H^*(S(X)/\On).\end{align} The images of the generators of the cohomology groups under these maps give Losik's characteristic classes.

\section{Expressions for some characteristic classes}\label{secexpr}

Now we describe the construction of the isomorphisms \eqref{isomH} following \cite{L10}. Recall the definition of the canonical Gelfand-Kazhdan form \cite{BR}. 
Let $U\subset\Real^n$ be an open subset. 
Let $\tau$ be a tangent
vector at $s\in S(U)$ and  let $s(u)$ be a curve on $S(U)$ such that
$\tau=\frac{ds}{du}(0)$. One can represent $s(u)$ by a smooth family
$k_u$ of germs at $0$ of regular at $0\in\Real^n$ maps $\Real^n\to U$,
i.e., $s(u)=j^\infty_0k_u$. Then  put
$$
\theta(\tau)=-j_0^\infty\frac{d}{du}(k_0^{-1}\circ k_u)(0).
$$
The form $\theta$ is canonical, i.e., it is invariant under the infinite order frame bundle maps induced by the diffeomorphisms of open subsets of $\Real^n$. Consequently, $\theta$ is a well-defined 1-form on $S(\pt_n)$ with values in $W_n$.
Now, each continuous cochain $c\in C^k(W_n)$ defines $\theta_c\in
\Omega^k(S(\pt))$ by the formula
\begin{equation}\label{mapforms}\theta_c(V_1,\dots, V_k)=c(\theta(V_1),\dots, \theta(V_k)),\end{equation}
where $V_1,\dots, V_k$ are vector fields on $S(\pt_n)$. The form $\theta$ satisfies the Maurer-Cartan equation, this implies that the map $c\mapsto \theta_c$ defines the homomorphism $C^*(W_n)\to \Omega^*(S(\pt_n))$ of the complexes. This homomorphism induces homomorphisms of the relative complexes.
All these cochain homomorphisms are isomorphisms that induce \eqref{isomH}.


Let us collect some information about Gelfand-Fuchs cohomology \cite{Fuchs} following \cite{L10}.
The complex
$C^*(W_n)$ is generated by the 1-forms $c^{i}_{j_1,\dots,j_r}$, $r=0,1,2\dots$, $1\leq i,j_1,\dots,j_r\leq n$,
where $$c^{i}_{j_1,\dots,j_r}(\xi)=\frac{\partial^r\xi^i}{\partial x^{j_1}\cdots \partial x^{j_r}}(0),\quad
\xi=\xi^i\partial_{x^i}\in W_n.$$ 
Let $$\gamma=\left(c^i_j\right), \quad \Psi^i_j=c^i_{jk}\wedge c^k, \quad \Psi=\left(\Psi^i_j\right).$$
The cochains $$\Psi_p=\text{tr}(\underbrace{\Psi\wedge\cdots\wedge\Psi}_{p \text{ times}}),\quad p=1,\dots,n,$$
are cocycles of $C^*(W_n,\GLn)$, and the cohomology classes of these cocycles generate $H^*(W_n,\GLn)$. The cohomology classes of $\Psi_p$ are called formal Chern classes.
 
There exist cochaines $\Gamma_1,\dots,\Gamma_n$ in $C^*(W_n)$ such that $d\Gamma_i=\Psi_i$. In particular, $\Gamma_1=c^i_i$. The cohomology $H^*(W_n)$ are isomorphic to the cohomology of the complex generated by  $\Gamma_1,\dots,\Gamma_n,\Psi_1,\dots,\Psi_n$. Finally, there exist cochains $\Lambda_{2k+1}$ ($1\leq 2k+1\leq n$) from $C^*(W_n,\On)$ such that $d\Lambda_{2k+1}=\Psi_{2k+1}$. The construction of these cochains may be found in \cite{L10}. The cohomology $H^*(W_n,\On)$ are isomorphic to the cohomology of the complex generated by the cochains $\Lambda_1,\Lambda_3,\dots,\Psi_1,\dots,\Psi_n$.

The equality $$\theta=\sum_{r\geq 0}\frac{1}{(j_1+\cdots+j_r)!}\theta^k_{j_1\dots j_r} x^{j_1}\cdots x^{j_r}\partial_{x^k}$$ defines on $S(\pt_n)$ canonical 1-forms $\theta^k_{j_1\cdots j_r}$. Applying \eqref{mapforms}, one gets  
$$ \theta^k_{j_1\dots j_r}=\theta_{c^k_{j_1\dots j_r}},$$ see \cite{L10}.
The explicit formulas for the canonical forms $\theta^k_{j_1\cdots j_r}$ for $r=0,1,2$ may be found in \cite{Kob61}.

Each  cohomology class from $H^*(W_n)$, $H^*(W_n,\GLn)$, and $H^*(W_n,\On)$ may be represented by cocycles that depend on the second jets of the vector fields \cite{Fuchs}. Together with the description of the form $\theta$ this shows that the images of the characteristic maps \eqref{chi}, \eqref{chi'}, \eqref{chi''} may be represented by forms on  the spaces of frames of the second order $S_2(X)$, $S_2(X)/\GLn$, and $S_2(X)/\On$, respectively. To simplify the exposition we will consider the characteristic classes with values in 
$H^*(S_2(X))$, $H^*(S_2(X)/\GLn)$, and $H^*(S_2(X)/\On)$.

Let $U,W\subset\Real^n$ be  open subsets such that $0\in W$. 
Let $f:W\to U$, $z^i=f^i(t^1,\dots,t^n)$, be a smooth map regular at the point $0$. The corresponding element of $S_2(U)$ has the
coordinates \begin{equation} \label{coordz} z^i_{j_1\dots j_r} =\frac{\partial^r f^i}{\partial t^{j_1}\cdots\partial t^{j_r}}(0),\quad r=0,1,2.
	\end{equation} 
	It is convenient to consider the following coordinate system:
	\begin{equation} \label{coordy} y^i=z^i,\quad y^i_{j}=z^i_{j},\quad y^i_{j k}=   v^p_jz^i_{pq}v^q_k,\end{equation}
	where $(v^i_j)$ is the inverse matrix to $(z^i_j)$.

	Suppose that we have
another open subset $V\subset\Real^n$ and a regular map $h:U\to V$.
Denote by $\alpha^k_{j_1\dots j_r}$, $r=0,1,2$, the coordinates on $S_2(V)$. Then we get the
induced map $\tilde h:S_2(U)\to S_2(V)$ given by
\begin{equation}\label{coordtrans1}\alpha^i=\alpha^i(y^1,\dots,y^n),\quad 
\alpha^i_j=\frac{\partial \alpha^i}{\partial y^p}y^p_j,\quad \alpha^i_{jk}=\frac{\partial y^p}{\partial \alpha^j }
\frac{\partial^2 \alpha^i }{\partial y^p \partial y^q}\frac{\partial y^q}{\partial \alpha^k}+
\frac{\partial y^p}{\partial \alpha^j }
\frac{\partial\alpha^i }{\partial y^s}\frac{\partial y^q}{\partial \alpha^k}y^s_{pq}.   \end{equation}

Using the results from \cite{Kob61}, we get
$$\theta^i_j=-v^i_kdy^k_j+y^i_{jk}dy^k.$$
Let $$\Gamma_1=\theta^i_i.$$
Note that $$d\Gamma_1=d\theta^i_i=dy^i_{ik}\wedge dy^k.$$
On the other hand, it is known \cite{Kob61} that 
$$d\theta^i_i=\theta^i_{ik}\wedge \theta^k.$$

The action of the group $\GLn$ on $S_2(U)$ with respect to the coordinates \eqref{coordy} is given by
\begin{equation}\label{action}
A:(y^i,y^i_j, y^i_{jk})\mapsto (y^i,y^i_pA^p_j, y^i_{jk}),\end{equation}
where $A=(A^i_j)\in\GLn$. 
This shows that the functions $y^i$, $y^i_{jk}$ may be considered as  coordinates on $S_2(U)/\GLn$.

\begin{defin}
	Let $X$ be a quotient. The first Chern-Losik class (CL class) is  the element of $H^2(S_2(X)/\GLn)$ defined by the image of the class  $[\Psi_1]\in H^2(W_n,\GLn)$ under the map~\eqref{chi'}.  
	\end{defin}

From the above we obtain

\begin{prop}
	The first CL class is the cohomology class of the form 
	$\theta^i_{ik}\wedge \theta^k$, which with respect to coordinates $y^i$, $y^i_{jk}$ is given by the form $dy^i_{ik}\wedge dy^k$.
\end{prop}

Consider the space $S_2(U)/\On$. Let $y^{ij}$ be the coordinates on the space of symmetric non-degenerate bilinear forms on $(\Real^n)^*$.  From \eqref{action} it follows that the functions $y^i$, $y^{ij}$, $y^i_{jk}$ may be considered as the coordinates on $S_2(U)/\On$. 
The projection $p:S_2(U)\to S_2(U)/\On$ is given by $$(y^i,y^i_j, y^i_{jk})\mapsto \left(y^i,\sum_l y^i_l y^k_l, y^i_{jk}\right).$$ 
A morphism $f:U\to V$ induces the map $\tilde f:S_2(U)/\On\to S_2(V)/\On$. Under $\tilde f$, the transformation of the coordinates $y^i$ and $y^i_{jk}$ is given by \eqref{coordtrans1}, and, moreover, it holds 
\begin{equation}\label{coordtrans2}
\alpha^{ij}=\frac{\partial \alpha^i}{\partial y^k}y^{kl} \frac{\partial \alpha^j}{\partial y^l}.
\end{equation}
Let $(y_{ij})$ denote the inverse matrix to $(y^{ij})$. Consider the following 1-form on $S(U)/\On$:
$$\Lambda_1=-\frac{1}{2}y_{ij}dy^{ij}+y^j_{jk}dy^k.$$ It is not hard to check 
\begin{lem} The 1-form $\Lambda_1$ is canonical, it satisfies $p^*\Lambda_1=\Gamma_1$, and $d\Lambda_1=d\Gamma_1$.
\end{lem}

Now we define two generalizations of the classical Godbillon-Vey class.

\begin{defin}
	Let $X$ be a quotient. The Godbillon-Vey-Losik class of $X$ with values in $H^{2n+1}(S_2(X))$ is the  cohomology class  $[\Gamma_1\wedge (d\Gamma_1)^n]\in H^{2n+1}(S_2(X))$.  
\end{defin}

\begin{defin}
	Let $X$ be a $\D_n$-space. The Godbillon-Vey-Losik class of $X$ with values in $H^{2n+1}(S_2(X)/\On)$ is the  cohomology class  $[\Lambda_1\wedge (d\Lambda_1)^n]\in H^{2n+1}(S_2(X)/\On)$.  
\end{defin}

It is clear that the projection $S_2(X)\to S_2(X)/\On$ preserves the GVL classes. 
From the above we obtain

\begin{prop}
	The GVL classes with values in $H^{2n+1}(S_2(X))$ is the cohomology class represented by the $(2n+1)$-form   that with respect to the coordinates $y^i$, $y^i_j$, $y^i_{jk}$ may be written as   
		$$-v^j_idy^i_j\wedge (dy^{k}_{k l}\wedge dy^l)^n.$$ 

\end{prop}

\begin{prop}
	The GVL classes with values in $H^{2n+1}(S_2(X)/\On)$ is the cohomology class represented by the $(2n+1)$-form that with respect to the coordinates $y^i$, $y^{ij}$, $y^i_{jk}$  may be written as 
	$$-\frac{1}{2}y_{ij}dy^{ij}\wedge (dy^{k}_{kl}\wedge dy^l)^n.$$

\end{prop}

\section{Characteristic classes associated to vector fields}

Let $M$ be a smooth manifold of dimension $n$. Let $\V$ be a set of vector fields on $M$. Denote by $P$  the pseudogroup of local diffeomorphisms of $M$ generated by the pseudogroups of the elements of $\V$.  We will denote the quotient $M/P$ simply by $M/\V$. We define the characteristic classes of  $\V$ to be the characteristic classes of $M/P$.

Let $U\subset\Real^n$ be an open subset. Let $V=V^i\partial_{y^i}$ be a vector field on $U$. The local flow of $V$ is extended to local flows on $S_2(U)$ and induces a vector field there. The same holds for the other bundles.  Define the functions $$\tilde V^i=V^i,\quad \tilde V^i_j=\frac{\partial V^i}{\partial y^p}y^p_j, \quad \tilde V^{ij}=\frac{\partial V^i}{\partial y^p}y^{pj}+\frac{\partial V^j}{\partial y^p}y^{ip}
\quad
\tilde V^i_{jk}=\frac{\partial^2 V^i}{\partial y^j\partial y^k}-
\frac{\partial V^p}{\partial y^j}y^i_{pk}-\frac{\partial V^p}{\partial y^k}y^i_{jp}+\frac{\partial V^i}{\partial y^p}y^p_{jk}.$$
Using \eqref{coordtrans1} it is easy to check that

\begin{prop}\label{proptildeV} The components of the extension $\tilde V$ of $V$ to $S_2(U)$ are $\tilde V^i$, $\tilde V^i_j$, 
	$\tilde V^i_{jk}$. 
		The components of the extension $\tilde V$ of $V$ to $S_2(U)/\GLn$ are $\tilde V^i$, 
$\tilde V^i_{jk}$. The components of the extension $\tilde V$ of $V$ to $S_2(U)/\On$ are $\tilde V^i$, $\tilde V^{ij}$, 
$\tilde V^i_{jk}$.  \end{prop}

The following obvious proposition shows how to check the (non-)triviality of the characteristic  classes. 

\begin{prop} Let $\V$ be a set of vector fields  on an $n$-dimensional smooth manifold $M$.  A cohomology class $[\alpha]\in H^k(S_2(M/\V))$  is trivial if and only if there exists a form
	$\omega\in\Omega^{k-1}(S_2(M))$ such that
	$$d\omega=\alpha$$ and
	$$L_{\tilde V}\omega=0$$ for all $V\in\V$. The same statements hold for the spaces $S_2(X)/\GLn$ and $S_2(X)/\On$.  \end{prop}

For the first CL class, the triviality may be checked using the following proposition.

\begin{prop}\label{proptrivcondCL} The first CL class is trivial in $H^2(S_2(M/\V)/\GLn)$ if and only if there is a cover $\{U_a\}$ of $M$ by coordinate neighborhoods, and
	functions $\eta_a$ on $S_2(U_a)/\GLn$ such that for each neighborhood $U_a$ it holds 
		\begin{equation}\label{propcond1}
	\frac{\partial V^i}{\partial y^i}+\tilde V\eta_a=c(a,V)=const, \quad V\in\V, \end{equation}
		and on each intersection $U_a\cap U_b$
		it holds 
		\begin{equation}\label{propcond2}d\eta_a=\frac{\partial y^p}{\partial \alpha^j }
		\frac{\partial^2 \alpha^j }{\partial y^p \partial y^q}dy^q  +d\eta_b,\end{equation}
where $y^i$ and $\alpha^i$ are coordinates on $U_a$ and $U_b$, respectively.	
	  \end{prop} 

{\bf Proof.} Suppose that the first CL class is trivial. By the previous proposition, there exists a 1-form $\omega$ on $S(M)/\GLn$ such that $L_{\tilde V}\omega=0$ for all $V\in\V$, and locally it holds $d\omega=dy^i_{ik}\wedge dy^k$. Suppose that the domains of the coordinate neighborhoods covering $M$ are simply connected. Then, for each coordinate neighborhood $U_a$, there exists a function $\eta_a$ on $S_2(U_a)/\GLn$ such that $$\omega=    y^i_{ik}dy^k+d\eta_a.$$
The equality \eqref{coordtrans1} and the assumption that  $\omega$ is a 1-form on $S_2(M)/\GLn$ imply \eqref{propcond2}. 
Let $V$ be a vector field on $U_a$. Using the above expression for $\tilde V$, it is easy to check that
$$L_{\tilde V}\omega=d\left(\frac{\partial V^i}{\partial y^i}+\tilde V\eta_a\right).$$
This implies \eqref{propcond1}.
Conversely, if the functions $\eta_a$ as above are given, then from \eqref{propcond2} and \eqref{coordtrans1} it follows that the forms $y^i_{ik}dy^k+d\eta_a$ define a 1-form $\omega$ on $S_2(M)/\GLn$, and  \eqref{propcond1} implies $L_{\tilde V}\omega=0$ for all $V\in\V$. This proves the proposition. \qed

Let $f$ be a (local) diffeomorphism of $M$. Suppose that $f$ may be included into a flow $f_t$, $f=f_1$. Let $V$ be the vector field of the flow $f_t$.

\begin{prop}\label{propw'} Let $G$ be one of the groups $\{e\}$, $\GLn$, $\On$.
	Let $\alpha\in\Omega^k (S_2(M)/G)$ be $\tilde f$-invariant. Then there exists a $\tilde f$-invariant form $\omega\in\Omega^{k-1}(S_2(M)/G)$ with $d\omega=\alpha$ if and only if there exists a form $\omega'\in\Omega^{k-1}(S_2(M)/G)$ such that $d\omega'=\alpha$ and $L_{\tilde V}\omega'=0$. 
\end{prop}
 {\bf Proof.} 
 It is enough to show that if  there
exists a $\tilde f$-invariant form $\omega\in\Omega^{k-1}(S_2(M)/G)$ with $d\omega=\alpha$, then there exists a form $\omega'\in\Omega^{k-1}(S_2(M)/G)$ such that $d\omega'=\alpha$ and $L_{\tilde V}\omega'=0$.
Let $\omega$ be  such a form. Consider the form 
\begin{equation}\label{Omega'}
\omega' = \int_{0}^{1}\tilde{f}^*_t(\omega)dt.
\end{equation}
 The form $\omega'$ satisfies $d\omega'=\alpha$, and it is $\tilde{f}_t$-invariant for each $t\in\mathbb{R}$, hence 
$L_{\tilde V}\omega'=0$. \qed

\begin{ex} Fix a positive real number $\alpha$. Consider  the  vector field $V_\alpha=V_\alpha(x)\frac{d}{dx}$ on $\Real$, where 
	\begin{equation}\label{Eq4.26}
	V_\alpha(x) = \left\{
	\begin{array}{cc}
	e^{-\frac{1}{|x|^\alpha}}, & \mbox{for } x \neq 0,\\
	0, & \mbox{for } x = 0.
	\end{array}
	\right.	
	\end{equation}
	
Let $X_\alpha=\Real/V_\alpha$. From the results of \cite{BGG} it follows that if $\alpha$ is an odd  integer, then the GVL class of $X_\alpha$ is non-trivial; if $\alpha$ is an even  integer, then the GVL class of $X_\alpha$ is trivial. This implies that the vector fields $V_{2m}$ and $V_{2k+1}$ cannot be conjugate by a diffeomorphism of $\Real$.
\end{ex}

\section{Simplified versions of some characteristic classes}\label{SimpVersOfClass}

The considered above characteristic classes of a quotient $X=M/P$ take values in the cohomology of the spaces $S_2(X)$, $S_2(X)/\GLn$, and  $S_2(X)/\On$. The dimensions of these spaces are big compared to the dimension of $X$. By that reason here we define the new versions of some of the above classes with values in certain new spaces associated to $X$.  
The above formulas for the GVL classes suggest first of all to consider new versions of these classes with values in the cohomology $S_2(X)/\SLn$  and $S_2(X)/(\SLnn)$, where $\SLnn$ is the subgroup of $\GLn$ consisting of matrices with the determinant $\pm 1$.  

For an open subset $U\subset\Real^n$,  we consider the local coordinates $y^i,y,y^i_{jk}$ on the space $S_2(U)/\SLn$ such the natural projection $S_2(U)\to S_2(U)/\SLn$ is given by
$$(y^i,y^i_j,y^i_{jk})\mapsto (y^i,\det(y^i_j),y^i_{jk}).$$
For a transformation $\alpha^i=\alpha^i(y^1,\dots,y^n)$, the coordinates $y^i$ and $y^i_{jk}$ are transformed according to 
\eqref{coordtrans1}. In addition it holds
$$\alpha=\det\left(\frac{\partial \alpha^i}{\partial y^j}\right)y.$$
It is clear that we may define the GVL  class with values in $H^3(S_2(X)/\SLn)$ in such a way that locally it is defined by the form
	$$-\frac{1}{y}dy\wedge (dy^{k}_{k l}\wedge dy^l)^n.$$

Next, consider the space $S_2(X)/(\SLnn)$ with the coordinates 
 $(y^i,x,y^i_{jk})$, $x=\ln|y|$. Then the GVL class with values in 
$H^3(S_2(X)/\SLn)$ 
is defined by the form
$$-dx\wedge (dy^{k}_{k l}\wedge dy^k)^n.$$

The sequence of projections 
$$S_2(X)\to S_2(X)/\On\to S_2(X)/\SLn \to S_2(X)/(\SLnn)$$
defines the sequence of the cohomology homomorphisms 
$$H^3(S_2(X)/(\SLnn))\to H^3(S_2(X)/\SLn)\to H^3(S_2(X)/\On)\to H^3(S_2(X))$$ that respects the different versions of the GVL classes.

Next, the formulas for the GVL classes suggest to consider the space $\tilde S_2(X)$ with the coordinates $y^i,y^i_j,y_i$, where $y_i=y^j_{ji}$. The coordinates $y_i$ are transformed by the rule
$$\alpha_k=\frac{\partial y^p}{\partial \alpha^j }
\frac{\partial^2 \alpha^j }{\partial y^p \partial y^q}\frac{\partial y^q}{\partial \alpha^k}+
\frac{\partial y^j}{\partial \alpha^k }
y_{j}.$$ We obtain the space $\tilde S_2(X)/\SLn$ with the coordinates $y^i,y,y_i$ and the GVL class with values in $H^3(\tilde S_2(X)/\SLn)$ is given by
$$-\frac{1}{y}dy\wedge (dy_i\wedge dy^i)^n,$$
i.e., in that case, GVL is defined by the volume form.
Likewise, consider  the space $\tilde S_2(X)/(\SLnn)$ with the coordinates $y^i,x,y_i$ and the GVL class given by 
$$-dx\wedge (dy_i\wedge dy^i)^n.$$ The space $\tilde S_2(X)/(\SLnn)$ is exactly the space $A(X)$ defined in the introduction.
Finally, we consider the space $B(X)=\tilde S_2(X)/\GLn$ with the coordinates $y^i,y_i$ and obtain the new version of the first CL class, which is now given by the symplectic form
$$dy_i\wedge dy^i.$$

Let $V$ be a vector field on an open subset $U\subset\Real^n$ with the components $V^i$. Define the following functions on $A(U)$:
$$\tilde V^i=V^i,\quad \hat V=\frac{\partial V^i}{\partial y^i}y, \quad  \tilde V_i = \dfrac{\partial^2V^j}{\partial y^j\partial y^i}-\dfrac{\partial V^j}{\partial y^i}y_j.$$ 

The following propositions are similar to Propositions~\ref{proptrivcondCL} and~\ref{proptildeV}, respectively.

\begin{prop}\label{proptildeVAB}  The components of the extension $\tilde V$ of $V$ to $A(U)$ are $\tilde V^i$, $\hat V$, 
	$\tilde V^i_{jk}$. 
	The components of the extension $\tilde V$ of $V$ to $B(U)$ are $\tilde V^i$, 
	$\tilde V_{i}$.   \end{prop}

\begin{prop}\label{proptrivcondCLB} Let $\V$ be a set of vector fields on a manifold $M$. The first CL class of $\V$ is trivial in $H^2(B(M/\V))$ if and only if there is a cover $\{U_a\}$ of $M$ by coordinate neighborhoods, and
	functions $\eta_a$ on $B(U_a)$ such that for each neighborhood $U_a$ it holds 
	\begin{equation}\label{propcond1B}
	\frac{\partial V^i}{\partial y^i}+\tilde V\eta_a=c(a,V)=const, \quad V\in\V, \end{equation}
	and on each intersection $U_a\cap U_b$
	it holds 
	\begin{equation}\label{propcond2B}d\eta_a=\frac{\partial y^p}{\partial \alpha^j }
	\frac{\partial^2 \alpha^j }{\partial y^p \partial y^q}dy^q  +d\eta_b,\end{equation}
	where $y^i$ and $\alpha^i$ are coordinates on $U_a$ and $U_b$, respectively.	
\end{prop}

\section{Characteristic classes of vector fields on 2-dimensional discs}\label{secex}

As examples we will consider  quotients of the open 2-dimensional disc $D^2\subset\Real^2$ by the flow  of a vector field $V$ on $D^2$. 

Let us fix a vector field $V$ on $D^2$. Denote by $V^1,V^2$ the components of $V$.
 Let $y^1,y^2$ be the coordinates on $D^2$ obtained as restriction of the standard coordinates on $\Real^2$.  Let $y^1,y^2,y_1,y_2$ be the coordinates on $B(D^2)=D^2\times\Real^2$. Let $\tilde V$ be the extension of $V$ to $B(D^2)$.  Proposition \ref{proptrivcondCLB} implies

\begin{prop} The first Chern-Losik class of a vector field $V$ on $D^2$ is trivial in $H^2(B(D^2/V))$ if and only if there exists a smooth function $G$ on 
	$D^2\times\Real^2$ satisfying the equation 
\begin{equation}\label{CLequation}
{\tilde V}G =  -\dfrac{\partial  V^1}{\partial y^1} - \dfrac{\partial  V^2}{\partial y^2} + R
\end{equation} for a constant $R\in\Real$.	
	
	\end{prop}

 Let us note that the function 
\begin{equation}\label{IntegralofCLequation}
	\dfrac{\partial V^i}{\partial y^i} - V^iy_i.
\end{equation}
is a first integral of equation \eqref{CLequation}.

In what follows we use the notation
$$r=r(y^1,y^2)=\sqrt{(y^1)^2+(y^2)^2}.$$
 
\subsection*{Example 1}

 Let  $f:[0,1)\to\Real$ be a smooth function non-vanishing on $(0,1)$. Recall that if $f$ admits an even  smooth extension to $(-1,1)$, then the function $f(r)$ on $D^2$ is smooth. Let us assume that $f$ satisfies this condition.  Then it holds $f^{(2n-1)}(0) = 0$ for all $n\in\mathbb{N}.$  Consider the vector field $V$ on $D^2$ with the components
\begin{equation}\label{VectField1}
			V^1(y^1,y^2) = f(r)y^1,\quad V^2(y^1,y^2) = f(r)y^2.
	\end{equation} 
\begin{theorem}\label{CLtheor1}
	The first CL class of the quotient $D^2/V$ is non-trivial in $H^2(B(D^2/V))$ 
	if and only if $f(0)=0$.
\end{theorem}
{\bf Proof.} 
The components of vector field $\tilde V$ are the following:
\begin{align*}
			\tilde{V}^1 &= V^1 = f(r)y^1,\\
		\tilde{V}^2 &= V^2 = f(r)y^2,\\
		\tilde{V}_1 &= \dfrac{f'(r)}{r}y^1(3-y_1y^1-y_2y^2) + f''(r)y^1 - f(r)y_1,\\
		\tilde{V}_2 &= \dfrac{f'(r)}{r}y^2(3-y_1y^1-y_2y^2) + f''(r)y^2 - f(r)y_2.
	\end{align*}
	The equation \eqref{CLequation} takes the form
	\begin{equation}\label{CLequation1}
	{\tilde V}G =  -2f(r)-rf'(r)+R.	\end{equation} 
	
	Fix a number $r_0\in (0,1)$. It is clear that the function
		$$G_0(r)=-\int\limits_{r_0}^{r}\dfrac{2f(\rho)+
		f'(\rho)\rho-R}{f(\rho)\rho}d\rho$$ is smooth on the open set $D^2\backslash\{0\}$ and it is easy the check that $G_0$ satisfies \eqref{CLequation1}.

	{\bf Case 1: $f(0)\neq 0$.} This means that $f(r)>0$ for all $r\in [0,1)$.  Let  $R=2f(0)$. We will prove that the function $G_0$ is smooth on $D^2$. 
		Let $r\in (-1,1)$ be considered as a variable.   We need to prove that the function $G_0(r)$, $r\in(-1,1)$ is smooth  and even.
			Since $f'(0)=0$, by Hadamard's Lemma, there exists a smooth function $g(r)$ such that $$f(r)=r^2g(r)+f(0).$$ Since $f$ is even, the function $g$ is even es well.
	We get that 
	$$G_0(r)=-\ln(f(r))+\ln(f(r_0))- \int\limits_{r_0}^{r}\dfrac{2\rho g(\rho)}{f(\rho)}d\rho.$$
	Since the function $\ln(f(r))$ is smooth and even on $(-1,1)$,  we only need to prove that the function
	$$h(r)=\int\limits_{r_0}^{r}\dfrac{2\rho g(\rho)}{f(\rho)}d\rho$$ is smooth and even.
	    It holds  
	    $$h'(r)=2\dfrac{1}{f(r)}r g(r),$$ 
	    which implies that $h(r)$ is smooth. It is obvious that the function $h(r)$ is even. 	Thus the function $G_0$ is smooth on $D^2$, i.e., it is a solution of equation  \eqref{CLequation1} smooth on $D^2\times \Real^2$.

	{\bf Case 2: $f(0) = 0$.} In this case we are going to prove that the equation \eqref{CLequation1} does not have any solution smooth on $D^2\times\Real^2$.
	
	\begin{lem}
		If $f(0)=0$, then the function $G_0$ has no smooth extension from  $D^2\backslash\{(0,0)\}$ to~$D^2$.
	\end{lem}
	{\bf Proof.}
	 Suppose that the function $G_0(r)$ has a smooth extension to $D^2$. Then the function $$\frac{2f(r)+f'(r)r-R}{f(r)r}$$ of the real variable $r\in(0,1)$ has a smooth extension to $[0,1)$. This implies $R=0$. Hence the function takes the form
	$$\dfrac{2}{r}+\dfrac{f'(r)}{f(r)}.$$ Since $f$ is non-vanishing on $(0,1)$, it is either positive or negative on $(0,1)$. Suppose for simplicity that $f$ is positive; the case of negative $f$ may be considered in the same way.
	Since $f(1/n) > 0$ for $n\in\mathbb{N}$, $n\geqslant 2$, by  Lagrange's mean value theorem, there exists $\xi_n\in (0,1/n)$ such that $$0<f(1/n) = f'(\xi_n)\cdot 1/n.$$ 
	Hence, $f'(\xi_n) >0$.  It is clear that $\xi_n\to +0$ as $n\to\infty$. Moreover, since $f'(\xi_n)>0$ and $f(\xi_n)>0$, it holds
	$$\dfrac{2}{\xi_n}+\dfrac{f'(\xi_n)}{f(\xi_n)}  > \dfrac{2}{\xi_n} \xrightarrow[n\to\infty]{} +\infty,$$
	i.e.,   $G_0$ has no smooth extension to $D^2$. \qed

Consider the vector field 
$$Y=y^2\frac{\partial}{\partial y^1}-y^1\frac{\partial}{\partial y_2}+y_2\frac{\partial}{\partial y_1}-y_1\frac{\partial}{\partial y_2}.$$	
	
\begin{lem} The equation \eqref{CLequation1} has a smooth solution on $D^2\times\Real^2$
if and only the system of equations 
\begin{equation}\label{CLequation2}{\tilde V}Q =  -2f(r)-rf'(r)+R,\quad YQ=0.\end{equation}
has a smooth solution on $D^2\times\Real^2$.
\end{lem}

{\bf Proof.} The flow $\varphi_\alpha$ of the vector field $Y$  consists of the transformations $$
\begin{pmatrix}
y^1\\ 
y^2
\end{pmatrix} \mapsto
\begin{pmatrix}
\cos\alpha & \sin\alpha\\
-\sin\alpha & \cos\alpha
\end{pmatrix}
\begin{pmatrix}
y^1 \\
y^2
\end{pmatrix}, \quad 
\begin{pmatrix}
y_1\\ 
y_2
\end{pmatrix} \mapsto
\begin{pmatrix}
\cos\alpha & \sin\alpha\\
-\sin\alpha & \cos\alpha
\end{pmatrix}
\begin{pmatrix}
y_1 \\
y_2
\end{pmatrix}.$$

Let $G$ be a solution of \eqref{CLequation1}. It is easy to check that $\varphi_\alpha^*$ preserves the equation \eqref{CLequation1}. Hence, for each $\alpha$, $\varphi_\alpha^*G$ is again a solution of \eqref{CLequation1}.
Consequently, the function $$Q=\frac{1}{2\pi}\int_0^{2\pi}(\varphi_\alpha^*G)d\alpha$$ is a solution of 
\eqref{CLequation2}. \qed

Let us consider the homogeneous system \begin{equation}\label{CLequation2hom}{\tilde V}Q = 0,\quad YQ=0. \end{equation} 

Consider the open set
$$U=(D^2\backslash \{(0,0)\})\times \Real^2.$$

\begin{lem}
	If Q is a function satisfying \eqref{CLequation2hom} and smooth on $U$, then it has a smooth extension to $D^2\times\Real^2$.
\end{lem}

{\bf Proof.}
 The vector fields $Y$ and $\tilde V$ are commuting and linearly independent at all points from the set $U$. Consequently, these vector fields define a foliation $\mathcal{T}$ on $U$. The solutions of \eqref{CLequation2hom} are exactly the functions constant along the leaves of the foliation $\mathcal{T}$.
It is easy to check that the functions
$$I_1=f(r)(2-y_1y^1-y_2y^2)+f'(r)r, \quad I_2=y_1y^2-y_2y^1$$
are solutions of the system \eqref{CLequation2hom} smooth on $D^2\times\Real^2$. 

Consider the map $$\psi:D^2\times\Real^2\to\Real^2,\quad \psi(y)=(I_1(y),I_2(y)).$$
For each $(y^1,y^2)\in D^2\backslash \{(0,0)\}$, let $$S_{(y^1,y^2)}=\{(y^1,y^2)\}\times\Real^2.$$ Since, for the corresponding $r$, $f(r)\neq 0$, the map
$$\psi|_{S_{(y^1,y^2)}}:S_{(y^1,y^2)}\to\Real^2$$ is a diffeomorphism.
Consequently, for each fixed $(c_1,c_2)\in\Real^2$, the map
$$\phi_{(c_1,c_2)}:D^2\backslash \{(0,0)\}\to (\psi|_U)^{-1}(c_1,c_2),\quad
\phi_{(c_1,c_2)}(y^1,y^2)=(\psi|_{S_{(y^1,y^2)}})^{-1}(c_1,c_2)$$
 is a diffeomorphism as well. This shows that each subset $(\psi|_U)^{-1}(c_1,c_2)\subset U$ is connected and consequently it coincides with a leaf of the foliation $\mathcal{T}$, i.e., the leaves of the foliation $\mathcal{T}$ are connected. Thus we see that each set  $S_{(y^1,y^2)}$ intersects each leaf of the foliation $\mathcal{T}$ at a single point. 
 
 Let $Q$ be a solution of \eqref{CLequation2hom} smooth on $U$.
 Fix a point $(y^1,y^2)\in D^2\backslash \{(0,0)\}$.
 Let $$H=\Big(\big(\psi|_{S_{(y^1,y^2)}}\big)^{-1}\Big)^*Q$$ be a smooth function on $\Real^2$. 
 	Then the equality \begin{equation}\label{eqQHI}
 	Q=H(I_1,I_2)
 	\end{equation}
 	holds on  	$S_{(y^1,y^2)}$. Consequently this equality holds  on each leaf of the foliation $\mathcal{T}$. Hence, the equality holds  on $U$. This implies that $H(I_1,I_2)$ is a smooth extension of $Q$ to $D^2\times \Real^2$. \qed

Suppose now that $G$ is a solution of \eqref{CLequation1} smooth on $D^2\times\Real^2$.
Then the function $$Q=G-G_0$$ is a solution of \eqref{CLequation2hom} smooth on $U$,
consequently, $Q$ has a smooth extension to $D^2\times\Real^2$. This implies that
the function $$G_0=G-Q,$$ has a smooth extension to $D^2\times\Real^2$, which gives a contradiction. \qed

\subsection*{Example 2}

Consider the vector field $V$ on $D^2$ with the following components:
$$V^1(y^1,y^2) = f(r)y^2,\quad
V^2(y^1,y^2) = -f(r)y^1,$$

where $f:(-1,1)\to\mathbb{R}$ is a smooth even function satisfying $f^{(2n-1)}(0) = 0$ for all $n\in\mathbb{N}$. 

\begin{theorem}
	The first CL class of the quotient $D^2/V$ is trivial in $H^2(B(D^2/V))$. \end{theorem}
{\bf Proof.}
The components of vector field $\tilde V$ generated by $V$ are
\begin{align*}
\tilde{V}^1 &= V^1 = f(r)y^2;\\
\tilde{V}^2 &= V^2 = -f(r)y^1;\\
\tilde{V}_1 &= f(r)y_2 + \dfrac{f'(r)}{r}y^1(y_2y^1-y_1y^2);\\
\tilde{V}_2 &= -f(r)y_1 + \dfrac{f'(r)}{r}y^2(y_2y^1-y_1y^2).
\end{align*} 
Hence the equation (\ref{CLequation})  takes the form
\begin{equation*}
\tilde VG= R,
\end{equation*}
where $R\in \Real$ is a constant. It is easy to check that the function $$G = y_2y^1-y_1y^2$$ is a solution of this system  with $R = 0$.
\qed

\vskip0.5cm


\end{document}